\input amstex
\documentstyle{amsppt}

\input label.def
\input degt.def

\def\paragraph{\subsection{}}
\copycounter\thm\subsection
\copycounter\equation\subsection
\newtheorem{Remark}\Remark\thm\endAmSdef
\problem\thm\endproclaim

{\catcode`\@11
\gdef\proclaimfont@{\sl}
\gdef\subsubheadfont@{\bf}
}

\def\dash{\item"\hfill--\hfill"}
\def\Dashes{\widestnumber\item{--}\roster}
\def\endDashes{\endroster}

\def\ie{\emph{i.e.}}
\def\eg{\emph{e.g.}}
\def\etc{\emph{etc}}
\def\cf.{\emph{cf\.}}
\def\via{\emph{via}}

\def\Cp#1{\Bbb P^{#1}}
\def\Rp#1{\Bbb P_\R^{#1}}
\def\XR{X_\R}
\def\discr{\operatorname{discr}}
\def\tsmash#1{\smash{\tilde#1}}
\def\tS{\tsmash{S}}
\def\tSigma#1{\tsmash{\Sigma^{#1}}}
\def\Ell{\bold{E}}
\let\CC=B
\def\tX{\tsmash{X}}
\def\tXR{\tsmash{X_\R}}
\def\tA#1{\tsmash{A_#1}}
\def\tB{\tsmash{B}}

\def\CL{\Cal L}
\def\CK{\Cal K}
\def\CM{\Cal M}
\def\CN{\Cal N}
\def\CS{\Cal S}
\let\onto\twoheadrightarrow
\let\into\hookrightarrow
\let\<\langle
\let\>\rangle

\let\Gf\varphi
\def\bA{\bold{A}}
\def\bD{\bold{D}}

\topmatter

\author
Alex Degtyarev
\endauthor

\address
Department of Mathematics,
Bilkent University,
06800 Ankara, Turkey
\endaddress

\email
degt\@fen.bilkent.edu.tr
\endemail

\title
Towards the generalized Shapiro and Shapiro conjecture
\endtitle

\dedicatory
To my teacher Oleg Viro at his 60th birthday
\enddedicatory

\abstract
We find a new, asymptotically better, bound $g\le\frac14d^2+O(d)$
on the genus of a curve that may violate the generalized total
reality conjecture. The bound covers all known cases except $g=0$
(the original conjecture).
\endabstract



\endtopmatter

\document

\section{Introduction}

The original (rational) total reality conjecture suggested by B.
and M.~Shapiro in 1993 states that, if all flattening points of a
regular curve $\Cp1\to\Cp{n}$ belong to the real line
$\Rp1\subset\Cp1$, then the curve can be made real by an
appropriate projective transformation of~$\Cp{n}$. (The
\emph{flattening points} are the points in the source~$\Cp1$ where
the first $n$ derivatives of the map are linearly dependent. In
the case $n=1$, a curve is a meromorphic function and the
flattening points are its critical points.) There is
a number of interesting and not always straightforward
restatements of this conjecture, in terms of the Wronsky map,
Schubert calculus, dynamical
systems, \etc. Although supported by extensive numerical evidence,
the conjecture proved extremely difficult to settle. It was not
before 2002 that the first result appeared, due to A.~Eremenko and
A.~Gabrielov~\cite{EG}, settling the case $n=1$, \ie,
meromorphic functions on~$\Cp1$. Later, a number of sporadic
results were announced, and the conjecture was proved in full
generality in 2005 by E.~Mukhin, V.~Tarasov, and
A.~Varchenko, see~\cite{MTV}.
The proof, revealing a deep
connection between Schubert calculus and theory of integrable
system, is based on the Bethe ansatz method in the Gaudin
model.

In the meanwhile, a number of generalizations of the conjecture
were suggested. In this paper, we deal with one of them,
see~\cite{ESS} and Problem~\ref{conj} below,
replacing
the source~$\Cp1$ with an arbitrary compact complex curve (but,
however, restricting~$n$ to~$1$, \ie, to the case of meromorphic
functions). Due to the lack of evidence, the authors chose to
state the assertion as a problem rather than a conjecture.

Recall that a \emph{real variety} is a complex
algebraic (analytic)
variety~$X$ supplied with a \emph{real structure},
\ie, an anti-ho\-lo\-morphic involution $c\:X\to X$. Given two
real varieties $(X,c)$ and $(Y,c')$, a regular map $f\:X\to Y$ is
called \emph{real} if it commutes with the real structures:
$f\circ c=c'\circ f$.

\problem[Problem \rm(see~\cite{ESS})]\label{conj}
Let $(C,c)$ be a real curve and let $f\:C\to\Cp1$ be
a regular
map such that\rom:
\roster
\item\local{conj.1}
all critical points and critical values of~$f$ are pairwise
distinct;
\item\local{conj.2}
all critical points of~$f$ are real.
\endroster
Is it true that $f$ is real with respect to an appropriate real
structure in~$\Cp1$?
\endproblem

The condition that the critical points of~$f$ are distinct
includes, in particular, the requirement that each critical point
is simple, \ie, has ramification index~$2$.

A pair of integers $g\ge0$, $d\ge1$ is said to have the
\emph{total reality property} if the answer to Problem~\ref{conj}
is in the affirmative for any curve~$C$ of genus~$g$ and map~$f$
of degree~$d$. At present, the total reality property is known for
the following pairs $(g,d)$:
\Dashes
\dash
$(0,d)$ for any $d\ge1$ (the original conjecture, see~\cite{EG});
\dash
$(g,d)$ for any $d\ge1$ and
$g>G_1(d):=\frac13(d^2-4d+3)$, see~\cite{ESS};
\dash
$(g,d)$ for any $g\ge0$ and $d\le4$, see~\cite{ESS}
and~\cite{DShapiro}.
\endDashes
The principle result of the present paper is the following
theorem.

\theorem\label{th.main}
Any pair $(g,d)$ with $d\ge1$ and $g$ satisfying the
inequality
$$
g>G_0(d):=\cases
k^2-2k,&\text{if $d=2k$ is even},\\
k^2-\dfrac{10}3k+\dfrac73,&\text{if $d=2k-1$ is odd}
\endcases
$$
has the total reality property.
\endtheorem

\Remark
Note that one has $G_0(d)-G_1(d)\le-\frac13(k-1)^2\le0$, where
$k=[\frac12(d+1)]$. Theorem~\ref{th.main} covers the values $d=2$,
$3$ and leaves only $g=0$ for $d=4$, reducing the generalized
conjecture to the classical one. The new bound is also
asymptotically better:
$G_0(d)=\frac14d^2+O(d)<G_1(d)=\frac13d^2+O(d)$.
\endRemark

\subsection{Contents of the paper}
In \S\ref{S.reduction}, we outline the reduction of
Problem~\ref{conj} to the question of existence of certain real
curves on the ellipsoid and restate Theorem~\ref{th.main} in the
new terms, see Theorem~\ref{th.curve}. In \S\ref{S.discr}, we
briefly
recall V.~V.~Nikulin's theory of discriminant forms and lattice
extensions. In \S\ref{S.Alexander}, we introduce a version of the
Alexander module of a plane curve suited to the study of the
resolution lattice in the homology of the double covering of the
plane ramified at the curve. Finally, in \S\ref{S.proof}, we prove
Theorem~\ref{th.curve} and, hence, Theorem~\ref{th.main}.

\subsection{Acknowledgements}
I am thankful to T.~Ekedahl and B.~Shapiro for the fruitful
discussions of the subject.
This work was completed during my
participation in the special semester on Real and Tropical
Algebraic Geometry held at \emph{Centre Interfacultaire
Bernoulli}, \emph{\'Ecole polytechnique f\'ed\'erale de Lausanne}.
I extend my gratitude to the organizers of the semester and to the
administration of \emph{CIB}.

\section{The reduction}\label{S.reduction}

We briefly recall the reduction of Problem~\ref{conj} to the
problem of existence of a certain real curve on the ellipsoid.
Details are found in~\cite{ESS}.

\paragraph
Denote by $\conj\:z\mapsto\bar z$
the standard real structure on $\Cp1=\C\cup\infty$.
The \emph{ellipsoid}~$\Ell$ is the quadric $\Cp1\times\Cp1$ with
the real structure $(z,w)\mapsto(\conj w,\conj z)$. (It is indeed
the real structure whose real part is homeomorphic to the
$2$-sphere.)

Let $(C,c)$ be a real curve and let $f\:C\to\Cp1$ be a holomorphic
map. Consider the \emph{conjugate map}
$\bar f=\conj\circ f\circ c\:C\to\Cp1$ and let
$$
\Phi=(f,\bar f)\:C\to\Ell.
$$
It is
straightforward that $\Phi$ is holomorphic and real (with respect
to the above real structure on~$\Ell$). Hence, the image $\Phi(C)$
is a real algebraic curve in~$\Ell$. (We exclude the possibility
that $\Phi(C)$ is a point as we assume $f\ne\const$,
\cf. Condition~\iref{conj}{conj.1}.) In particular, the image
$\Phi(C)$ has bi-degree $(d',d')$ for some $d'\ge1$.

\lemma[Lemma \rm(see~\cite{ESS})]
A holomorphic map $f\:C\to\Cp1$ is real with respect to some real
structure on~$\Cp1$ if and only if there is a M\"obius
transformation $\Gf\:\Cp1\to\Cp1$ such that $\bar f=\Gf\circ f$.
\qed
\endlemma

\corollary[Corollary \rm(see~\cite{ESS})]\label{(1,1)}
A holomorphic map $f\:C\to\Cp1$ is real with respect to some real
structure on~$\Cp1$ if and only if the image $\Phi(C)\subset\Ell$,
see above,
is a curve of bi-degree $(1,1)$.
\qed
\endcorollary

\paragraph
Let $p\:\Ell\to\Cp1$ be the projection to the first factor.
In general, the map $\Phi$ as above splits into a ramified
covering~$\alpha$ and a generically one-to-one map~$\beta$,
$$
\Phi\:C@>\alpha>>C'@>\beta>>\Ell,
$$
so that $d=\deg f=d'\deg\alpha$, where $d'=\deg(p\circ\beta)$ or,
alternatively, $(d',d')$ is the bi-degree of the image
$\Phi(C)=\beta(C')$. Then, $f$ itself splits into~$\alpha$ and
$p\circ\beta$. Hence, the critical values of~$f$ are those
of~$p\circ\beta$ and the images under~$p\circ\beta$ of the
ramification points of~$\alpha$. Thus, if $f$ satisfies
Condition~\iref{conj}{conj.1}, the splitting cannot be proper,
\ie, either $d=\deg\alpha$ and $d'=1$ or $\deg\alpha=1$ and
$d=d'$. In the former case, $f$ is real with respect to some real
structure on~$\Cp1$, see Corollary~\ref{(1,1)}. In the latter
case, assuming that the critical points of~$f$ are real,
Condition~\iref{conj}{conj.2}, the image $B=\Phi(C)$ is a curve of
genus~$g$ with $2g+2d-2$ real ordinary cusps (type~$\bA_2$
singular points, the images of the critical points of~$f$) and all
other singularities with smooth branches.

Conversely, let $B\subset\Ell$
be a real curve of bi-degree $(d,d)$, $d>1$, and genus~$g$
with $2g+2d-2$ real ordinary cusps and all other singularities
with smooth branches, and let $\rho\:\tsmash{B}\to B$ be the
normalization of~$B$. Then $f=p\circ\rho\:\tsmash{B}\to\Cp1$ is a
map that satisfies Conditions~\iref{conj}{conj.1}
and~\ditto{conj.2} but is not real with respect to any real
structure on~$\Cp1$; hence, the pair $(g,d)$ does not have the
total reality property.

As a consequence, we obtain the following statement.

\theorem[Theorem \rm(see~\cite{ESS})]\label{th.reduction}
A pair $(g,d)$ has the total reality property if and only if there
does \emph{not} exist a real curve $B\subset\Ell$ of degree~$d$
and genus~$g$ with $2g+2d-2$ real ordinary cusps and all other
singularities with smooth branches.
\qed
\endtheorem

Thus, Theorem~\ref{th.main} is equivalent to the following
statement, which is actually proved in the paper.

\theorem\label{th.curve}
Let $\Ell$ be the ellipsoid, and let $\CC\subset\Ell$ be a
real curve of bi-degree $(d,d)$ and genus~$g$ with $c=2d+2g-2$
real ordinary cusps and other singularities
with smooth branches. Then
$g\le G_0(d)$, see Theorem~\ref{th.main}.
\endtheorem

\Remark
It is worth mentioning that the bound $g>G_1(d)$ mentioned in the
introduction is purely complex: it is derived from the adjunction
formula for the virtual genus of a curve $B\subset\Ell$ as in
Theorem~\ref{th.reduction}. On the contrary, the
proof of the conjecture for the case $(g,d)=(1,4)$ found
in~\cite{DShapiro} makes an essential use of the real structure,
as an elliptic curve with eight ordinary cusps in $\Cp1\times\Cp1$
does exist! Our proof of Theorem~\ref{th.curve} also uses the
assumption that all cusps are real.
\endRemark

\paragraph\label{s.ph}
In general, a curve~$B$ as in Theorem~\ref{th.curve} may have
rather complicated singularities. However, as the proof below is
essentially topological, we follow S.~Yu.~Orevkov \cite{Orevkov}
and
perturb~$B$ to a
real \emph{pseudo-holomorphic} curve
with ordinary nodes (type $\bA_1$) and ordinary cusps
(type~$\bA_2$) only. By the genus formula, the number of nodes of
such a curve is
$$
n=(d-1)^2-g-c=d^2-4d-1-3g.
\eqtag\label{eq.n}
$$

\section{Discriminant forms}\label{S.discr}

In this section, we cite the techniques and a few results of
Nikulin~\cite{Nikulin}. Most proofs are found
in~\cite{Nikulin}; they are omitted.

\paragraph\label{s.lattice}
A \emph{lattice} is a finitely generated free abelian group~$L$
equipped with a symmetric bilinear form $b\:L\otimes L\to\ZZ$. We
abbreviate $b(x,y)=x\cdot y$ and $b(x,x)=x^2$.
As the transition matrix between two integral bases
has determinant $\pm1$, the
determinant $\det L\in\ZZ$
(\ie, the determinant
of the Gram matrix of~$b$ in
any
basis of~$L$)
is well defined.
A lattice~$L$ is called
\emph{nondegenerate} if $\det L\ne0$; it is called
\emph{unimodular} if $\det L=\pm1$ and \emph{$p$-unimodular} if
$\det L$ is prime to~$p$ (where $p$ is a prime).

To fix the notation, we use $\sigma_+(L)$, $\sigma_-(L)$, and
$\sigma(L)=\sigma_+(L)-\sigma_-(L)$ for, respectively, the
positive and negative inertia indices and the signature of a
lattice~$L$.

\paragraph\label{s.discr}
Given a lattice~$L$,
the bilinear form extends to $L\otimes\Q$. If
$L$ is nondegenerate, the dual group $L^*=\Hom(L,\ZZ)$ can
be regarded as the subgroup
$$
\bigl\{x\in L\otimes\Q\bigm|
 \text{$x\cdot y\in\ZZ$ for all $x\in L$}\bigr\}.
$$
In particular, $L\subset L^*$ and the quotient $L^*/L$
is a finite group; it is called the {\it discriminant group\/}
of~$L$ and is denoted by $\discr L$ or~$\CL$. The
group~$\CL$
inherits from $L\otimes\Q$ a symmetric bilinear form
$\CL\otimes\CL\to\Q/\ZZ$,
called the {\it discriminant form};
when
speaking about the discriminant groups, their
(anti-)isomorphisms, \etc., we always assume that the discriminant
form
is taken
into account.
The following properties are straightforward:
\roster
\item\local{CL.1}
the discriminant form is nondegenerate, \ie, the associated
homomorphism $\CL\to\Hom(\CL,\Q/Z)$ is an isomorphism;
\item\local{CL.2}
one has $\#\CL=\mathopen|\det L\mathclose|$;
\item\local{CF.3}
in
particular, $\CL=0$ if and only if $L$ is unimodular.
\endroster

Following Nikulin, we denote by $\ell(\CL)$ the minimal number of
generators of a finite abelian group~$\CL$. For a prime~$p$, we
denote by~$\CL_p$ the $p$-primary part of~$\CL$ and let
$\ell_p(\CL)=\ell(\CL_p)$. Clearly, for a lattice~$L$ one has
\roster
\item[4]\local{CL.4}
$\rank L\ge\ell(\CL)\ge\ell_p(\CL)$ (for any prime~$p$);
\item\local{CL.5}
$L$ is $p$-unimodular if and only if $\CL_p=0$.
\endroster

\paragraph\label{s.extension}
An \emph{extension} of a lattice~$S$ is another lattice~$M$
containing~$L$. All lattices below are assumed nondegenerate.

Let $M\supset S$ be a finite index extension of a
lattice~$S$. Since $M$ is also a
lattice, one has monomorphisms $S\into M\into M^*\into S^*$. Hence,
the quotient $\CK=M/S$ can be regarded as a subgroup of
the discriminant $\CS=\discr S$; it
is called the \emph{kernel} of the extension $M\supset S$. The
kernel is an isotropic subgroup, \ie, $\CK^\perp\subset\CK$, and
one has $\CM=\CK^\perp\!/\CK$. In particular, in view
of~\iref{s.discr}{CL.1}, for any prime~$p$ one has
$$
\ell_p(\CM)\ge\ell_p(\CL)-2\ell_p(\CK).
$$

Now, assume that $M\supset S$ is a \emph{primitive} extension,
\ie, the quotient $M/S$ is torsion free. Then the
construction above applies to the finite index extension
$M\supset S\oplus N$, where
$N=S^\perp$, giving rise to the kernel
$\CK\subset\CS\oplus\CN$. Since both~$S$ and~$N$ are primitive
in~$M$, one has $\CK\cap\CS=\CK\cap\CN=0$; hence, $\CK$ is the
graph of an anti-isometry~$\kappa$ between certain subgroups
$\CS'\subset\CS$ and $\CN'\subset\CN$. If $M$ is unimodular, then
$\CS'=\CS$ and $\CN'=\CN$, \ie, $\kappa$
is an anti-isometry $\CS\to\CN$.
Similarly, if $M$ is $p$-unimodular for
a certain
prime~$p$, then $\CS'_p=\CS_p$ and $\CN'_p=\CN_p$, \ie, $\kappa$
is an anti-isometry $\CS_p\to\CN_p$. In particular,
$\ell(\CS)=\ell(\CN)$ (respectively, $\ell_p(\CS)=\ell_p(\CN)$).
Combining these observations
with~\iref{s.discr}{CL.4},
we arrive at the following statement.

\lemma\label{rk.Sperp}
Let $p$ be a prime, and let $L\supset S$ be a $p$-unimodular
extension of a nondegenerate lattice~$S$. Denote by~$\tS$ the
primitive hull of~$S$ in~$L$, and let~$\CK$ be the kernel of the
finite index extension $\tS\supset S$. Then
$\rank S^\perp\ge\ell_p(\CS)-2\ell_p(\CK)$.
\qed
\endlemma

\section{The Alexander module}\label{S.Alexander}

Here, we discuss (a version of) the Alexander module of a plane
curve and its relation to the resolution lattice in the homology
of the double covering of the plane ramified at the curve.


\paragraph\label{s.A}
Let $\pi$ be a group, and let $\kappa\:\pi\onto\Z_2$ be an
epimorphism. Denote $K=\Ker\kappa$ and define
the \emph{Alexander module} of~$\pi$ (more precisely,
of~$\kappa$) as the $\Z[\Z_2]$-module $A_\pi=K/[K,K]$, the
generator $t$ of~$\Z_2$ acting \via\
$x\mapsto[\bar t^{-1}\bar x\bar t]\in A_\pi$, where $\bar t\in\pi$
and $\bar x\in K$ are some representatives of~$t$ and~$x$,
respectively. (We simplify the usual definition and consider only
the case needed in the sequel. A more general version and further
details can be found in A.~Libgober~\cite{Libgober}.)

Let $B\subset\Cp1\times\Cp1$ be an irreducible curve of even
bi-degree $(d,d)=(2k,2k)$, and let
$\pi=\pi_1(\Cp1\times\Cp1\sminus B)$. Recall that
$\pi/[\pi,\pi]=\Z_{2k}$; hence, there is a unique epimorphism
$\kappa\:\pi\onto\Z_2$. The resulting Alexander module
$A_B=A_\pi$ will be
called the \emph{Alexander module} of~$B$. The \emph{reduced
Alexander module}~$\tA{B}$ is the kernel of the canonical
homomorphism $A_B\to\Z_k\subset\pi/[\pi,\pi]$. There is a natural exact
sequence
$$
0@>>>\tA{B}@>>>A_B@>>>\Z_k@>>>0
\eqtag\label{eq.A}
$$
of $\Z[\Z_2]$-modules (where the $\Z_2$-action on $\Z_d$ is
trivial).
The following statement is essentially contained in
O.~Zariski~\cite{Zariski}.

\lemma\label{splitting}
The exact sequence~\eqref{eq.A} splits\rom: one has
$A_B=\tA{B}\oplus\Ker(1-t)$,
where $t$ is the generator of~$\Z_2$.
Furthermore,
$\tA{B}$ is a finite group free of $2$-torsion, and the action
of~$t$ on $\tA{B}$ is \via\ the multiplication by $(-1)$.
\endlemma

\proof
Since $A_B$ is a finitely generated abelian group,
to prove that it is finite and free of $2$-torsion
it suffices to show that
$\Hom_\Z(\tA{B},\Z_2)=0$. Assume the contrary. Then the
$\Z_2$-action in the $2$-group $\Hom_\Z(\tA{B},\Z_2)$ has a fixed
non-zero element, \ie, there is an equivariant epimorphism
$\tA{B}\onto\Z_2$. Hence, $\pi$ factors to a group~$G$ that is an
extension $0\to\Z_2\to G\to\Z_{2k}\to0$. The group~$G$ is
necessarily abelian and it is strictly larger than
$\Z_{2k}=\pi/[\pi,\pi]$. This is a contradiction.

Since $\tA{B}$ is finite and free of $2$-torsion, one can divide
by~$2$ and there is a splitting $\tA{B}=\tA{{}^+}\oplus\tA{{}^-}$,
where $\tA{{}^\pm}=\Ker[(1\pm t)\:\tA{B}\to\tA{B}]$. Then, $\pi$
factors to a group~$G$ that is a central extension
$0\to\tA{{}^+}\to G\to\Z_{2k}\to0$, and as above one concludes
that $\tA{{}^+}=0$, \ie, $t$ acts on~$\tA{B}$ \via\ $(-1)$.

Pick a representative $a'\in A_B$ of a generator of
$\Z_k=A_B/\tA{B}$. Then, obviously,
$(1-t)a'\in\tA{B}$, and replacing $a'$
with $a'+\frac12(1-t)a'$, one obtains a $t$-invariant
representative $a\in\Ker(1-t)$. The multiple $ka\in\tA{B}$ is both
invariant and skew-invariant; since $\tA{B}$ is free of
$2$-torsion, $ka=0$ and the sequence splits.
\endproof

\paragraph\label{s.covering}
Let $B\subset\Cp1\times\Cp1$ be an irreducible curve of even
bi-degree $(d,d)=(2k,2k)$ and with simple singularities only.
Consider the double covering $X\to\Cp1\times\Cp1$ and denote
by~$\tX$ the minimal resolution of singularities of~$X$. Let
$\tB\subset\tX$ be the proper pull-back of~$B$, and let
$E\subset\tX$ be the exceptional divisor contracted by the
blow-down $\tX\to X$.

Recall that the minimal resolution of a simple surface singularity
is diffeomorphic to its perturbation, see, \eg,~\cite{Durfee}.
Hence, $\tX$ is diffeomorphic to the double covering of
$\Cp1\times\Cp1$ ramified at a nonsingular curve. In particular,
$\pi_1(\tX)=0$ and one has
$$
b_2(X)=\chi(X)-2=8k^2-8k+6,
\quad
\sigma(X)=-4k^2.
\eqtag\label{eq.sigma}
$$

\paragraph\label{s.Sigma}
Denote $L=H_2(\tX)$. We regard~$L$ as a lattice \via\ the
intersection index pairing on~$\tX$. (Since $\tX$ is simply
connected, $L$ is a free abelian group. It is a unimodular lattice
due to the Poincar\'e duality.)
Let $\Sigma\subset L$ be the sublattice spanned
by the components of~$E$, and let $\tSigma{}\subset L$ be the
primitive hull of~$\Sigma$. Recall that $\Sigma$ is a negative
definite lattice.
Let, further, $h_1,h_2\subset L$ be the classes of the
pull-backs of a pair of generic generatrices of $\Cp1\times\Cp1$,
so that $h_1^2=h_2^2=0$, $h_1\cdot h_2=2$.

\lemma\label{CK}
If a curve
$B$ as above is irreducible, then there are natural isomorphisms
$\tA{B}=\Hom_\Z(\CK,\Q/\Z)=\Ext_\Z(\CK,\Z)$, where $\CK$ is the kernel
of the
extension $\tSigma{}\supset\Sigma$.
\endlemma

\proof
One has $A_B=H_1(\tX\sminus(\tB+E))$ as a group, the $\Z_2$-action
being induced by the deck translation of the covering. Hence, by
the Poincar\'e--Lefschetz duality,
$A_B$ is the cokernel of the inclusion homomorphism
$i^*\:H^2(\tX)\to H^2(\tB+E)$.

On the other hand, there is an orthogonal (with respect to the
intersection index form in~$\tX$) decomposition
$H_2(\tB+E)=\Sigma\oplus\<b\>$, where $b=k(h_1+h_2)$ is the class
realized by the divisorial pull-back of~$B$ in~$\tX$. The cokernel
of the restriction $i^*\:H^2(X)\to\<b\>^*$ is a cyclic group
$\Z_k$ fixed by the deck translation. Hence, in view of
Lemma~\ref{splitting},
$$
\tA{B}=\Coker[i^*\:H^2(\tX)\to H^2(E)]=\Coker[L^*\to\Sigma^*]=
 \discr\Sigma/\CK^\perp.
$$
(We use the splitting $L^*\onto\tSigma*\to\Sigma^*$, the first map
being an epimorphism as $L/\tSigma{}$ is torsion free.)
Since the discriminant form is nondegenerate,
see~\iref{s.discr}{CL.1},
one has
$\discr\Sigma/\CK^\perp=\Hom_\Z(\CK,\Q/\Z)$.
Since $\CK$ is a finite group, applying the functor
$\Hom_\Z(\CK,\,\cdot\,)$ to the short exact sequence
$0\to\Z\to\Q\to\Q/\Z\to0$, one obtains an isomorphism
$\Hom_\Z(\CK,\Q/\Z)=\Ext_\Z(\CK,\Z)$.
\endproof

\corollary\label{abelian}
In the notation of Lemma~\ref{CK}, if $B$ is irreducible and the
group
$\pi_1(\Cp1\times\Cp1\sminus B)$ is abelian, then $\CK=0$.
\qed
\endcorollary

\corollary\label{CK<=}
In the notation of Lemma~\ref{CK}, if $B$ is an irreducible curve
of bi-degree $(d,d)$, $d=2k\ge2$, then $\CK$ is free of $2$-torsion
and
$\ell(\CK)\le d-2$.
\endcorollary

\proof
Due to Lemma~\ref{CK}, one can replace~$\CK$ with $\tA{B}$. Then,
the statement on the $2$-torsion is given by
Lemma~\ref{splitting}, and it suffices to estimate the numbers
$\ell_p(\tA{B})=\ell(\tA{B}\otimes\Z_p)$ for odd primes~$p$.

Due to the Zariski--van Kampen theorem~\cite{vanKampen} applied to
one of the two rulings of $\Cp1\times\Cp1$, there is an
epimorphism
$\pi_1(L\sminus B)=F_{d-1}\onto\pi_1(\Cp1\times\Cp1\sminus B)$,
where $L$ is a generic generatrix of $\Cp1\times\Cp1$ and
$F_{d-1}$ is the free group on $d-1$ generators. Hence, $A_B$ is a
quotient of the Alexander module
$$
\tsize
A_{F_{d-1}}=\Z[\Z_2]/(t-1)\oplus\bigoplus_{d-2}\Z[\Z_2].
$$
For an odd prime~$p$, there is a splitting
$A_{F_{d-1}}\otimes\Z_p=A_p^+\oplus A_p^-$ (over the field~$\Z_p$)
into the eigenspaces of the action of~$\Z_2$,
and, due to Lemma~\ref{splitting},
the group $\tA{B}\otimes\Z_p$ is a quotient of
$A_p^-=\bigoplus_{d-2}\Z_p$.
\endproof

\Remark\label{rem.ph}
All statements in this section hold for pseudo-holomorphic curves
as well, \cf.~\ref{s.ph}.
For Corollary~\ref{CK<=}, it suffices to assume that $B$
is a small perturbation of an algebraic curve of bi-degree
$(d,d)$. Then, one still has an epimorphism
$F_{d-1}\onto\pi_1(\Cp1\times\Cp1\sminus B)$, and the proof
applies literally.
\endRemark

\section{Proof of Theorem~\ref{th.main}}\label{S.proof}

As explained in~\S\ref{S.reduction}, it suffices to prove
Theorem~\ref{th.curve}. We consider the cases of $d$ even
and $d$ odd separately.

\paragraph\label{s.c}
Let $B\subset\Cp1\times\Cp1$ be an irreducible curve of even
bi-degree $(d,d)$, $d=2k$. Assume that all singularities
of~$B$ are simple and let $\tX$ be the minimal resolution
of singularities of the double covering $X\to\Cp1\times\Cp1$
ramified at~$B$, \cf.~\ref{s.covering}. As in \ref{s.Sigma},
consider the unimodular lattice $L=H_2(\tX)$.

Let $c\:\tX\to\tX$ be a real structure on~$\tX$, and denote
by~$L^\pm$ the $(\pm1)$-eigenlattices of the induced
involution~$c_*$ of~$L$. The following statements are well known:
\roster
\item\local{c.1}
$L^\pm$ are the orthogonal complements of each other;
\item\local{c.2}
$L^\pm$ are $p$-unimodular for any odd prime~$p$;
\item\local{c.3}
one has $\sigma_+(L^+)=\sigma_+(L_-)-1$.
\endroster
Since also
$\sigma_+(L^+)+\sigma_+(L^-)=\sigma_+(L)=2k^2-4k+3$,
see~\eqref{eq.sigma}, one arrives at
$\sigma_+(L^+)=\sigma_+(L^-)-1=(k-1)^2$ and, further, at
$$
\rank L^-=(7k^2-6k+5)-\sigma_-(L^+).
\eqtag\label{eq.rank}
$$

\Remark
The common proof of Property~\iref{s.c}{c.3} uses the Hodge
structure. However, there is another (also very well known) proof
that also applies to almost complex manifolds. Let $\tXR=\Fix c$
be the real part of~$\tX$. Then, the normal bundle of $\tXR$
in~$\tX$ is $i$ times its tangent bundle; hence, the normal Euler
number $\tXR\circ\tXR$ equals $(-1)$ times the index of any
tangent vector field on~$\tXR$, \ie, $-\chi(\tXR)$. Now, one has
$\sigma(L^+)-\sigma(L^-)=\tXR\circ\tXR=-\chi(\tXR)$
(by the Hirzebruch $G$-signature theorem) and
$\rank L^+-\rank L^-=\chi(\XR)-2$ (by the Lefschetz
fixed point theorem). Adding the two equations,
one obtains~\iref{s.c}{c.3}.
\endRemark

\subsection{The case of $d=2k$ even}\label{s.2k}
Perturbing, if necessary, $B$ in the class of real
pseudo-holomorphic curves, see~\ref{s.ph}, one can assume that
all singularities of~$B$ are
$c$ real
ordinary cusps and
$n$
ordinary nodes,
where
$$
c=2d+2g-2\quad\text{and}\quad
n=d^2-4d-1-3g,
\eqtag\label{eq.c,n}
$$
see Theorem~\ref{th.reduction} and~\eqref{eq.n}.
Let $n=r+2s$, where $r$ and~$s$
are the numbers of, respectively, real nodes and pairs of
conjugate nodes.

\paragraph\label{s.sing.pts}
Consider the double covering~$\tX$, see~\ref{s.covering}, lift
the real structure on~$\Ell$ to a real structure~$c$ on~$\tX$, and
let $L^\pm\subset L$ be the corresponding eigenlattices,
see~\ref{s.c}. In the notation of~\ref{s.Sigma}, let
$\Sigma^\pm=\Sigma\cap L^\pm$.
Then
\Dashes
\dash
each real cusp of~$B$ contributes a sublattice~$\bA_2$
to~$\Sigma^-$;
\dash
each real node of~$B$ contributes a sublattice $\bA_1=[-2]$
to~$\Sigma^-$;
\dash
each pair of conjugate nodes contributes $[-4]$ to~$\Sigma^-$ and
$[-4]$ to~$\Sigma^+$.
\endDashes
In addition, the classes~$h_1$,
$h_2$ of two generic generatrices of~$\Ell$
span a hyperbolic plane orthogonal to~$\Sigma$, see~\ref{s.Sigma}.
It contributes
\Dashes
\dash
a sublattice $[4]\subset L^-$ spanned by $h_1+h_2$, and
\dash
a sublattice $[-4]\subset L^+$ spanned by $h_1-h_2$.
\endDashes
(Recall that any real structure reverses the
canonical complex orientation of
pseudo-holomorphic curves.)

\paragraph\label{s.estimates}
All sublattices of $L^+$ described above are negative definite;
hence, their total rank $s+1$ contributes to $\sigma_-(L^+)$. The
total rank $2c+r+s+1$ of the sublattices of~$L^-$ contributes to
the rank of $S^-=\Sigma^-\oplus[4]\subset L^-$. Due
to~\eqref{eq.rank}, one has
$$
2c+n+2+\rank S^\perp\le7k^2-6k+5,
\eqtag\label{eq.2k}
$$
where $S^\perp$ is the orthogonal complement of~$S^-$ in~$L^-$.
All summands of~$S^-$ other than $\bA_2$ are $3$-unimodular, whereas
$\discr\bA_2$ is the group $\Z_3$ spanned by an element of square
$\frac13\bmod\Z$. Let $\tsmash{S^-}\supset S^-$ and
$\tSigma{}\supset\Sigma$ be the primitive hulls, and denote
by~$\CK^-$ and~$\CK$ the kernels of the corresponding finite index
extensions, see~\ref{s.extension}. Clearly,
$\ell_3(\CK^-)\le\ell_3(\CK)$ and, due to Corollary~\ref{CK<=}
(see also Remark~\ref{rem.ph}), one has
$\ell_3(\CK)\le d-2$. Then, using Lemma~\ref{rk.Sperp},
one
obtains $\rank S^\perp\ge c-2(d-2)$, and combining the last
inequality with~\eqref{eq.2k}, one arrives at
$$
3c+n-2(d-2)\le7k^2-6k+3.
$$
It remains to substitute the expressions for~$c$ and~$n$
given by~\eqref{eq.c,n} and solve for~$g$ to get
$$
g\le k^2-2k+\frac23.
$$
Since $g$ is an integer, the last inequality implies
$g\le G_0(2k)$ as in Theorem~\ref{th.curve}.

\subsection{The case of $d=2k-1$ odd}\label{s.2k-1}
As above, one can assume that $B$ has $c$ real ordinary cusps and
$n=r+2s$ ordinary nodes, see~\eqref{eq.c,n}. Furthermore, one can
assume that $c>0$, as otherwise $g=0$ and $d=1$. Then, $B$ has a
real cusp and, hence, a real smooth point~$P$.
Let $L_1$, $L_2$ be the two generatrices of~$\Ell$
passing through~$P$. Choose~$P$ generic, so that each~$L_i$,
$i=1,2$, intersects~$B$ transversally at $d$ points, and consider
the real curve $B'=B+L_1+L_2$ of even bi-degree $(2k,2k)$,
applying to it the same double covering arguments as above.
In
addition to the nodes and cusps of~$B$, the new curve~$B'$ has
$(d-1)$ pairs of conjugate nodes and a real triple (type~$\bD_4$) point
at~$P$ (with one real and two complex conjugate branches). Hence,
in addition to the classes listed in~\ref{s.sing.pts}, there are
\Dashes
\dash
$(d-1)$ copies of $[-4]$ in each~$\Sigma^+$, $\Sigma^-$ (from the
new conjugate nodes),
\dash
a sublattice
$[-4]\subset\Sigma^+$ (from the type~$\bD_4$ point), and
\dash
a sublattice $\bA_3\subset\Sigma^-$
(from the type~$\bD_4$ point).
\endDashes
Thus, inequality~\eqref{eq.2k}
turns into
$$
2c+n+2(d-1)+4+2+\rank S^\perp\le7k^2-6k+5.
$$
We will show that $\rank S^\perp\ge c$. Then, substituting the
expressions for~$c$ and~$n$, see \eqref{eq.c,n},
and
solving the resulting inequality in~$g$, one would obtain
$g\le G_0(2k-1)$, as required.

\paragraph
In view of Lemma~\ref{rk.Sperp}, in order to prove that
$\rank S^\perp\ge c$, it suffices to show that $\ell_3(\CK)=0$
(\cf. similar arguments in~\ref{s.estimates}).

Perturb~$B'$ to
a pseudo-holomorphic curve~$B''$, keeping the cusps of~$B'$ and
resolving the other singularities. (It
would suffice to resolve the singular points resulting from the
intersection $B\cap L_1$.) Then, applying the Zariski--van Kampen
theorem~\cite{vanKampen} to the ruling containing~$L_1$, it is easy
to show that the fundamental
group $\pi_1(\Cp1\times\Cp1\sminus B'')$ is cyclic.
Indeed,
let~$U$ be a small tubular neighborhood of~$L_1$ in
$\Cp1\times\Cp1$, and let $L''\subset U$ be a generatrix
transversal to~$B''$. Obviously, the epimorphism
$\pi_1(L_1''\sminus B'')\onto\pi_1(\Cp1\times\Cp1\sminus B'')$
given by the Zariski--van Kampen theorem factors through
$\pi_1(U\sminus B'')$, and the latter group is cyclic.

On the
other hand, the new double covering $\tsmash{X''}\to\Cp1\times\Cp1$
ramified
at~$B''$ is diffeomorphic to~$\tX$, and the diffeomorphism can be
chosen identical over the union of a collection of Milnor balls
about the
cusps of~$B'$. Thus, since $\discr\bA_1$ and $\discr\bD_4$ are
$2$-torsion groups, the perturbation does not change
$\CK\otimes\Z_3$, and Corollary~\ref{abelian}
(see also Remark~\ref{rem.ph}) imply that $\CK\otimes\Z_3=0$.
\qed

\refstyle{C}


\Refs

\ref{DE}
\by A.~Degtyarev, T.~Ekedahl, I.~Itenberg, B.~Shapiro, M.~Shapiro
\paper On total reality of meromorphic functions
\jour Ann. Inst. Fourier
\vol 57
\yr 2007
\issue 5
\pages 2015--2030
\endref\label{DShapiro}

\ref{Du}
\by A.~H.~Durfee
\paper Fifteen characterizations of rational double points
 and simple critical points
\jour Enseign. Math. (2)
\vol 25
\yr 1979
\issue 1--2
\pages 131--163
\endref\label{Durfee}

\ref{ESS}
\by T.~Ekedahl, B.~Shapiro, M.~Shapiro
\paper First steps towards total reality of meromorphic functions
\jour Moscow Math. J.
\toappear
\endref\label{ESS}

\ref{EG}
\by A.~Eremenko, A.~Gabrielov
\paper Rational functions with real critical points and the B. and M.
Shapiro conjecture in real enumerative geometry
\jour Ann. of Math. (2)
\vol 155
\yr 2002
\issue 1
\pages 105--129
\endref\label{EG}

\ref{vK}
\by E.~R.~van~Kampen
\paper On the fundamental group of an algebraic curve
\jour  Amer. J. Math.
\vol   55
\yr    1933
\pages 255--260
\endref\label{vanKampen}

\ref{MTV}
\by E.~Mukhin, V.~Tarasov, A.~Varchenko
\paper The B. and M. Shapiro conjecture in real algebraic
geometry and the Bethe ansatz
\finalinfo\tt arXiv:math.AG/0512299
\endref\label{MTV}

\ref{L}
\by A.~Libgober
\paper
Alexander modules of plane algebraic curves
\jour
Contemporary Math.
\vol    20
\yr 1983
\pages  231--247
\endref\label{Libgober}

\ref{N}
\by V.~V.~Nikulin
\paper Integral symmetric bilinear forms and some of their
geometric applications,
\jour Izv. Akad. Nauk SSSR Ser. Mat.
\vol 43
\yr 1979
\pages 111--177
\lang Russian
\transl\nofrills English transl. in
\jour Math. USSR--Izv.
\vol 43
\yr 1980
\pages 103--167
\endref\label{Nikulin}

\ref{O}
\by S.~Yu.~Orevkov
\paper Classification of flexible $M$-curves of degree $8$ up to isotopy
\jour GAFA 
\vol 12
\yr 2002
\issue 4
\pages 723--755
\endref\label{Orevkov}


\ref{Z}
\by O.~Zariski
\paper On the irregularity of cyclic multiple planes
\jour Ann. Math.
\vol    32
\yr     1931
\pages  485--511
\endref\label{Zariski}

\endRefs
\enddocument